\documentclass {article}
\usepackage{authblk}
\usepackage{etex}
\usepackage[utf8]{inputenc}
\usepackage[english]{babel}
\usepackage{amsmath}
\usepackage{amsthm}
\usepackage{amssymb}
\usepackage{sectsty}
\usepackage{titlesec}
\usepackage{color}
\usepackage[color,matrix,arrow]{xy}
\usepackage{amsgen}
\usepackage{amstext}
\usepackage{amsbsy}
\usepackage{amsopn}
\usepackage{amsfonts}
\usepackage{eepic}
\usepackage{graphicx}
\usepackage{epsf}
\usepackage{pstricks}
\usepackage{enumerate}
\usepackage{eqnarray}
\usepackage{faktor}
\xyoption{all}

\usepackage{pgf}
\usepackage{tikz}
\usetikzlibrary{automata, arrows.meta, positioning,decorations.pathmorphing}

\newcommand{\footrecall}[1]{%
} 
\usetikzlibrary{snakes,shapes,arrows,automata}

\titleformat*{\section}{\large\bfseries}
\titleformat*{\subsection}{\normalsize \bfseries}

\newcommand{\N}{\mathbb{N}}
\newcommand{\Z}{\mathbb{Z}}

\newcommand{\Orb}{\text{Orb}}
\newcommand{\Aut}{\text{Aut}}

\newcommand{\Rat}{\text{Rat}}
\newcommand{\Alg}{\text{Alg}}
\newcommand{\CF}{\text{CF}}
\newcommand{\Rec}{\text{Rec}}

\newcommand{\mc}{\mathcal}

\theoremstyle{definition}
\newtheorem{theorem}{Theorem}[section]
\newtheorem{corollary}[theorem]{Corollary}

\newtheorem{proposition}[theorem]{Proposition}
\newtheorem{problem}[theorem]{Problem}
\newtheorem{lemma}[theorem]{Lemma}
\newtheorem{example}[theorem]{Example}
\newtheorem{remark}[theorem]{Remark}
 
\begin{document}
 
% --------------------------------------------------------------
%                         Começar aqui
% --------------------------------------------------------------
 
\title{Algebraic and context-free subsets of subgroups}

\author{Andr\'e Carvalho \thanks{andrecruzcarvalho@gmail.com}}
\affil{Centre of Mathematics, University of Porto

Rua do Campo Alegre, s/n

4169-007 Porto, 
Portugal}

\date{}

\maketitle
\begin{abstract}
We study the relation between the structure of algebraic and context-free subsets of a group $G$ and that of a finite index subgroup $H$. Using these results, we prove that a kind of Fatou property, previously studied by Berstel and Sakarovitch in the context of rational subsets and by Herbst in the context of algebraic subsets, holds for context-free subsets if and only if the group is virtually free. We also exhibit a counterexample to a question of Herbst concerning this property for algebraic subsets.
\end{abstract}

\section{Introduction}
Rational and recognizable subsets of groups are natural generalizations of finitely generated and finite index subgroups, respectively. Over the years, they have been studied from different points of view. From the structure viewpoint, presumably the most important result is Benois' Theorem, which provides us with a description of rational subsets of free groups in terms of reduced words. A description of rational subsets of free-abelian groups as the semilinear sets of $\Z^m$ is also known (see \cite{[Ber79]} for references). Deep connections with the algebraic structure of the group have also been proved in many contexts. For example, rational (resp. recognizable) subgroups have been proved to be exactly the finitely generated (resp. finite index) ones, and some classes of groups can be described through the classification of their subsets (see for example \cite[Theorem 3.1]{[Her91]}, where several descriptions of virtually cyclic groups are given and Theorem \ref{vfree charact} below, where virtually free groups are characterized).  Another application of the study of subsets of groups is the generalization of some decision problems concerning finitely generated subgroups, such as the membership problem, the intersection problem or the generalized conjugacy problem (see, for example \cite{[KSS07],[LS08],[LS11]}).

In   \cite{[BS86]}, the authors notice that, using the techniques from \cite{[AS75]}, one can deduce that, given a group $G$ and a subgroup $H\leq G$, then 
\begin{align}
\label{fatou rat}
\Rat(H)=\{K\subseteq H\mid K\in \Rat(G)\},
\end{align}
where $\Rat(G)$ denotes the class of rational subsets of the group $G$. They call property (\ref{fatou rat}) \emph{a kind of Fatou property for groups.} 
Notice that it is clear that $\Rat(H)\subseteq \{K\subseteq H\mid K\in \Rat(G)\}$. The difficult part is to prove the reverse inclusion.

While more complex, the context-free counterparts of rational and recognizable subsets, respectively algebraic and context-free subsets, also yield interesting results.
In \cite{[Her91]}, Herbst studied these subsets and was able to characterize groups for which context-free subsets coincide with rational subsets as virtually cyclic groups. In the same paper, Herbst also proved the property (\ref{fatou rat}) for algebraic subsets in case $H$ is a finite index normal subgroup of $G$ and posed the question of whether this would hold in general. Later, Herbst proved that that was the case if $G$ is a virtually free group. We will exhibit a counterexample to this question. 
We will also consider the same question for recognizable and context-free subsets and prove that, in the first case, the property holds if and only if $H$ is a finite index subgroup of $G$ and in the latter it holds for all $H\leq_{f.g} G$ if and only if $G$ is virtually free. To achieve this, we prove some results relating the structure of algebraic and context free-subsets of a group $G$ with the structure of the corresponding subsets of a finite index subgroup $H$, obtaining structural results similar to the ones obtained for rational and recognizable subsets in \cite{[Gru99],[Sil02b]}.

The paper is organized as follows: in Section \ref{sec: prelim}, we present basic definitions and results concerning subsets of subgroups that will be used throughout the paper. In Section \ref{sec: finite index}, we prove a structural result relating context-free and algebraic subsets of a group $G$ and those of a subgroup of finite index $H$ and prove that some decidability questions on $G$ are equivalent to the corresponding questions on $H$.
Finally, in Section \ref{sec: fatou}, we will consider the property (\ref{fatou rat}) for recognizable, context-free and algebraic subsets proving that, for recognizable subsets it holds if and only if $H$ has finite index in $G$; for context-free subsets it holds  for all $H\leq_{f.g} G$ if and only if $G$ is virtually free; and for algebraic subsets it does not hold in general, answering a question of Herbst \cite{[Her91],[Her92]}.

\section{Preliminaries}
\label{sec: prelim}
We will now present the basic definitions and results on rational, algebraic and context-free subsets of groups. For more detail, the reader is referred to \cite{[Ber79]} and \cite{[BS21]}.

The set $\{1,\ldots, n\}$ will be denoted by $[n]$.
Let $G=\langle A\rangle$ be a finitely generated group, $A$ be a finite generating set, $\tilde A=A\cup A^{-1}$ and $\pi:\tilde A^*\to G$ be the canonical (surjective) homomorphism. 

A subset $K\subseteq G$ is said to be \emph{rational} if there is some rational language $L\subseteq \tilde A^*$ such that $L\pi=K$ and \emph{recognizable} if $K\pi^{-1}$ is rational. 
Given a word in $w\in \tilde A^*$ and a letter $a\in \tilde A$,  we denote by $n_a(w)$ the number of occurrences of $a$  in $w$.

We will denote by $\Rat(G)$ and $\Rec(G)$ the class of rational and recognizable subsets of $G$, respectively. Rational subsets generalize the notion of finitely generated subgroups.

\begin{theorem}[\cite{[Ber79]}, Theorem III.2.7]
\label{AnisimovSeifert}
Let $H$ be a subgroup of a group $G$. Then $H\in \Rat(G)$ if and only if $H$ is finitely generated.
\end{theorem}

Similarly, recognizable subsets generalize the notion of finite index subgroups.

\begin{proposition}
\label{rec fi}
Let $H$ be a subgroup of a group $G$. Then $H\in \Rec(G)$ if and only if $H$ has finite index in $G$.
\end{proposition}

In fact, if $G$ is a group and $K$ is a subset of $G$ then $K$ is recognizable if and only if $K$ is a (finite) union of cosets of a subgroup of finite index.

In case the group $G$ is a free group with basis $A$ with surjective homomorphism $\pi:\tilde A^*\to G$, we define the set of reduced words of $L\subseteq \tilde A^*$ by $$Red(L)=\{w\in \tilde A^*\mid w \text{ is reduced and there exists $u\in L$ such that $u\pi=w\pi$} \}.$$ 
Benois' Theorem provides us with a useful characterization of rational subsets in terms of reduced words representing the elements in the subset.
\begin{theorem}[Benois]
Let $F$ be a finitely generated free group with basis $A$. Then, a subset of $Red(\tilde A^*)$ is a rational language of $\tilde A^*$ if and only if it is a rational subset of $F$.
\end{theorem}

The structure of rational and recognizable subsets of a certain group can be described by that of a finite index subgroup, as proved independently by Grunschlag and Silva.
\begin{proposition}\cite{[Gru99],[Sil02b]}
\label{gru-silva rational}
Let $G$ be a finitely generated group and $H\leq_{f.i.} G$. If $G$ is the disjoint union $G=\cup_{i=1}^n Hb_i$, then $Rat(G)$ consists of all subsets of the form 
\begin{align*}
\bigcup_{i=1}^{n} L_ib_i \qquad (L_i\in Rat(H))
\end{align*}
\end{proposition}

\begin{proposition}\cite{[Gru99],[Sil02b]}
\label{gru-silva recognizable}
Let $G$ be a finitely generated group and $H\leq_{f.i.} G$. If $G$ is the disjoint union $G=\cup_{i=1}^n Hb_i$, then $Rec(G)$ consists of all subsets of the form 
\begin{align*}
\bigcup_{i=1}^{n} L_ib_i \qquad (L_i\in Rec(H))
\end{align*}
\end{proposition}

A natural generalization of these concepts concerns the class of context-free languages. 
A subset $K\subseteq G$ is said to be \emph{algebraic} if there is some context-free language $L\subseteq \tilde A^*$ such that $L\pi=K$ and \emph{context-free} if $K\pi^{-1}$ is context-free. 
We will denote by $\Alg(G)$ and $\CF(G)$ the class of algebraic and context-free subsets of $G$, respectively.  It follows from \cite[Lemma 2.1]{[Her91]} that $\CF(G)$ and $\Alg(G)$ do not depend on the alphabet $A$ or the surjective homomorphism $\pi$.

It is obvious from the definitions that $\Rec(G)$, $\Rat(G)$, $\CF(G)$ and $\Alg(G)$ are closed under union, since both rational and context-free languages are closed under union. The intersection case is distinct: from the fact that rational languages are closed under intersection, it follows that $\Rec(G)$ must be closed under intersection too. However $\Rat(G)$, $\Alg(G)$ and $\CF(G)$ might not be. Another important closure property is given by the following lemma from \cite{[Her91]}.

\begin{lemma}
\label{cfletra}
\cite[Lemma 4.1]{[Her91]} Let $G$ be a finitely generated group, $R\in \Rat(G)$ and $C \in \{\Rat, \Rec, \Alg, \CF\}$.
If $K \in C(G)$, then $KR, RK \in C(G).$
\end{lemma}
This lemma will be used often in this paper in the particular case where $R$ is a singleton.

For a finitely generated group $G$, it is immediate from the definitions that $$\Rec(G)\subseteq \CF(G) \subseteq \Alg(G)$$ and that $$\Rec(G)\subseteq \Rat(G) \subseteq \Alg(G).$$ It is proved in \cite{[Her91]} that 
$$\CF(G)=\Alg(G) \Leftrightarrow \CF(G)=\Rat(G) \Leftrightarrow \text{ G is virtually cyclic.}$$

However, there is no general inclusion between $\Rat(G)$ and $\CF(G)$. For example, if $G$ is virtually abelian, then $\CF(G)\subseteq \Alg(G)= \Rat(G)$ (and the inclusion is strict if the group is not virtually cyclic) and if the group is virtually free, then $\Rat(G)\subseteq \CF(G)$ (see \cite[Lemma 4.2]{[Her91]}).

In the case of the free group $F_n$ of rank $n\geq 1$, Herbst proves in  \cite[Lemma 4.6]{[Her91]} an analogue of Benois' Theorem for context-free subsets, proving that for a subset $K\subseteq F_n$, then $K\in \CF(F_n)$ if and only if the set of reduced words representing elements of $K$ is context-free.

It follows from the definition of a context-free subsets and the fact that we can decide membership of a word in a context-free language that, if we are given as input a context-free grammar generating $K\pi^{-1}$ for some subset $K$ and a word $w$ representing an element $g\in G$, we can decide if $g\in K$.
% in case $K\pi^{-1}$ is computable for some context-free subset $K$, then we can decide the membership problem in a context-free subset of a group $G$.
 
 %%%%%%%%%%%%%%%%%%%%%%%%%%%%%%%%%%%%%%%%%%%%%%%%%%%%%%%%%%%%%%%%%%%%%%%%%%%%%%%%%%%%%
 \section{Finite index subgroups}
 \label{sec: finite index}
 In this section we will study how the structure of $\CF(G)$ and $\Alg(G)$ is related to the structure of $\CF(H)$ and $\Alg(H)$, where $H$ is a finite index subgroup of $G$. The structural results obtained are similar to the ones concerning $\Rat(G)$ and $\Rec(G)$ (Propositions \ref{gru-silva rational} and \ref{gru-silva recognizable}).

%We will start by dealing with context-free subsets. 
We will start by stating some technical lemmas that will be useful throughout the paper. 

\begin{lemma} \label{inters cf rec}
Let $G$ be a finitely generated group and $K_1\in \CF(G)$ (resp. $K_1\in \Alg(G)$) and $K_2\in \Rec(G)$. Then $K_1\cap K_2\in \CF(G)$ (resp. $K_1\cap K_2 \in \Alg(G)$).
\end{lemma}
\noindent\textit{Proof.} Let $G=\langle A \rangle$ and $\pi:A^*\to G$ be a surjective homomorphism. If $K_1\in \CF(G)$, then $K_1\pi^{-1}$ is context-free and so is $(K_1\cap K_2)\pi^{-1}=K_1\pi^{-1}\cap K_2\pi^{-1}$, since context-free languages are closed under intersection with a rational language. 

Also, if $K_1\in \Alg(G)$ and $L$ is a context-free language such that $L\pi=K_1$, then $L\cap K_2\pi^{-1}$ is context-free and $(L\cap K_2\pi^{-1})\pi=K_1\cap K_2$. Hence, $K_1\cap K_2\in \Alg(G)$.
\qed\\

Let $G$ be a finitely generated group, $H\leq_{f.g.} G$ and $K\subseteq G$. Corollary 4.4 in \cite{[Her91]} states that if $K\in CF(G)$, then $K\cap H\in CF(H)$. The next lemma follows immediately from this.

\begin{lemma} Let $G$ be a finitely generated group and $H\leq_{f.g.} G$. Then 
\begin{align*}
\{K\subseteq H\mid K\in \CF(G)\}\subseteq \CF(H).
\end{align*}
\label{cf partial}
\end{lemma}

We will now prove the same result for algebraic subsets in the case $H$ has finite index in $G$. It follows from \cite[Proposition 5.4]{[Her91]} that if $H\trianglelefteq_{f.i.} G$ is a normal subgroup of a group $G$ and $K\subseteq H$  belongs to $\Alg(G)$, then it must belong to $\Alg(H)$. We will now show that the normality hypothesis may be removed.
\begin{lemma} \label{fi desce}
Let $G$ be a finitely generated group and $H\leq_{f.i}G$. 
$$\{K\subseteq H\mid K\in \Alg(G)\}\subseteq \Alg(H).$$
\end{lemma}
\noindent\textit{Proof.} Let  $K\in \Alg(G)$ be such that $K\subseteq H$. There exists a normal subgroup $F\leq H$ such that $F\trianglelefteq_{f.i.} G$ (and so $F\trianglelefteq_{f.i.} H$). Then $H$ has a decomposition as a disjoint union
 $$H= Fb_1\cup\cdots \cup Fb_n,$$
 for some $b_i\in H$ and $K$ can be written as a disjoint union of the form
 $$K=K\cap H=\bigcup_{i=1}^n Fb_i\cap K.$$
 We will prove that for every $i\in [n]$, $Fb_i\cap K\in \Alg(H)$, which suffices since $\Alg(H)$ is closed under union. 

So, let $i\in[n]$ and write $K_i=Fb_i\cap K$. Then $K_ib_i^{-1}\subseteq F$. Since $F$ has finite index in $G$, then $F\in \Rec(G)$, and so, by Lemma \ref{cfletra}, $Fb_i\in \Rec(G)$. Since $K\in \Alg(G)$, then by Lemma \ref{inters cf rec}, it follows that $K_i\in \Alg(G)$, which implies that $K_ib_i^{-1}\in \Alg(G)$. Since $K_ib_i^{-1}\subseteq F$ and $F$ is a finite index normal subgroup of $G$, then  \cite[Proposition 5.4]{[Her91]} yields that  $K_ib_i^{-1}\in \Alg(F)$. Hence  $K_ib_i^{-1}\in \Alg(H)$ and $K_i\in \Alg(H)$.  
\qed\\

Notice that, similarly to what happens in the rational case, the reverse inclusion holds for every finitely generated subgroup $H$ (not necessarily of finite index). 

The next lemma is similar to \cite[Proposition 5.5(a)]{[Her91]}. The only difference is that we remove the hypothesis of normality.
\begin{lemma} \label{fi sobe}
Let $G$ be a finitely generated group and $H\leq_{f.i}G$. Then $\CF(H)\subseteq \CF(G)$.
\end{lemma}
\noindent\textit{Proof.} Proceeding as in the proof of Lemma \ref{fi desce}, we write
 $$K=K\cap H=\bigcup_{i=1}^n Fb_i\cap K,$$
 where $F\trianglelefteq_{f.i} H$.
 We will prove that for every $i\in [n]$, $Fb_i\cap K\in \CF(G)$, which suffices since $\CF(G)$ is closed under union. 

So, let $i\in[n]$ and write $K_i=Fb_i\cap K$. Then $K_ib_i^{-1}\subseteq F\leq H$. Since $F$ has finite index in $H$, then $F\in \Rec(H)$, and so, by Lemma \ref{cfletra}, $Fb_i\in \Rec(H)$. Since $K\in \CF(H)$, then by Lemma \ref{inters cf rec}, it follows that $K_i\in \CF(H)$, which, again by  Lemma \ref{cfletra}, implies that $K_ib_i^{-1}\in \CF(H)$. By Lemma \ref{cf partial}, we have that  $K_ib_i^{-1}\in \CF(F)$. Using  \cite[Proposition 5.5(a)]{[Her91]}, we obtain that  $K_ib_i^{-1}\in \CF(G)$, which means that  $K_i\in \CF(G)$, by  Lemma \ref{cfletra}.
\qed\\

\begin{corollary}
Let $G$ be a finitely generated group and $H\leq_{f.i}G$. Then $$\{K\subseteq H\mid K\in \CF(G)\}= \CF(H).$$
\end{corollary}

We are now able to prove a structural result analogous to Propositions \ref{gru-silva rational} and \ref{gru-silva recognizable}, but for context-free subsets. This gives an explicit description of $\CF(G)$ based on $\CF(H)$, where $H$ is a finite index subgroup of $G$.
\begin{proposition}
\label{cf finite extension}
Let $G$ be a finitely generated group and $H\leq_{f.i.} G$. If $G$ is the disjoint union $G=\cup_{i=1}^n Hb_i$, then $\CF(G)$ consists of all subsets of the form 
\begin{align*}
\bigcup_{i=1}^{n} L_ib_i \qquad (L_i\in \CF(H))
\end{align*}
\end{proposition}
\noindent\textit{Proof.} Let $K\subseteq G$ be a subset of the form $\cup_{i=1}^{n} L_ib_i$ with $L_i\in \CF(H)$. By Lemma \ref{fi sobe}, we have that $L_i\in \CF(G)$ for all $i\in [n]$ and so $L_ib_i\in \CF(G)$, by Lemma \ref{cfletra}. Since $\CF(G)$ is closed under union, then $K\in \CF(G)$. 

Conversely, let $K\in \CF(G)$. Then $K$ can be written as a disjoint union $$K=K\cap G=\bigcup_{i=1}^n (Hb_i\cap K).$$
Put $K_i=Hb_i\cap K$ and let $i\in [n]$.
Since $H\leq_{f.i.} G$, then by Proposition \ref{rec fi}, we have that $H\in \Rec(G)$, and so, $Hb_i\in \Rec(G)$, by Lemma \ref{cfletra}. It follows from  Lemma \ref{inters cf rec} that $K_i\in \CF(G)$ and again by  Lemma \ref{cfletra}, $L_i=K_ib_i^{-1}\in \CF(G)$ and $L_i\subseteq H$. By Lemma \ref{cf partial}, we deduce that $L_i\in \CF(H)$ and thus
$$K=K\cap G=\bigcup_{i=1}^n K_i=\bigcup_{i=1}^n L_ib_i.$$
\qed\\

 Proceeding in the same manner, we can obtain the same result for algebraic sets. 
\begin{proposition}
\label{alg finite extension}
Let $G$ be a finitely generated group and $H\leq_{f.i.} G$. If $G$ is the disjoint union $G=\cup_{i=1}^n Hb_i$, then $\Alg(G)$ consists of all subsets of the form 
\begin{align*}
\bigcup_{i=1}^{n} L_ib_i \qquad (L_i\in \Alg(H))
\end{align*}
\end{proposition}

We remark that, combining \cite[Lemma 4.6]{[Her91]} with Proposition \ref{cf finite extension}, we get a description of context-free subsets of a finitely generated virtually free group.
\begin{corollary}
Let $G$ be a finitely generated virtually free group, $F=F_A$ be a finite index free subgroup of $G$ with a free basis $A$, and $G=\cup_{i=1}^n Fb_i$ be a decomposition of $G$ as a disjoint union of cosets of $F$.
Then $\CF(G)$ consists of the sets of the form
$$\bigcup_{i=1}^n L_ib_i,$$
where $L_i\subseteq \tilde A^*$ is such that  $Red(L_i)$ is a context-free language of $\tilde A^*$.
\end{corollary}

The proof of the first item of the following corollary is essentially the same as the proof of \cite[Lemma 4.4]{[Sil02b]}.  
\begin{corollary}
\label{c1}
Let $G$ be a group and $H\leq_{f.i.}G$. Then we have that 
\begin{enumerate}[i.]
\item $\CF(G)$ (resp. $\Alg(G)$) is closed under intersection if and only if $\CF(H)$  (resp. $\Alg(H)$) is closed under intersection;
\item $\CF(G)$  (resp. $\Alg(G)$) is closed under complement if and only if $\CF(H)$   (resp. $\Alg(H)$) is closed under complement;
\end{enumerate}
\end{corollary}
\noindent\textit{Proof.} We will only present the proof of the context-free case, as the algebraic case is entirely analogous.

Write $G$ as a disjoint union $G=\cup_{i=1}^n Hb_i$.
We start by proving i. 

Suppose that $\CF(G)$ is closed under intersection and let $K,K'\in \CF(H)$. Then, by Lemma \ref{fi sobe}, $K,K'\in \CF(G)$, and so $K\cap K'\in \CF(G)$. Since $K\cap K'\subseteq H$, by Lemma \ref{cf partial}, $K\cap K'\in \CF(H)$. 

Conversely, suppose that $\CF(H)$ is closed under intersection and take $K,K'\in \CF(G)$. Then 
\begin{align*}
K=\bigcup_{i=1}^n L_ib_i \quad\text{ and }\quad K'=\bigcup_{i=1}^n L_i'b_i,
\end{align*} 
for some $L_i, L_i'\in \CF(H)$, $i\in[n]$. Since the cosets $Hb_i$ are disjoint, it follows that $$K\cap K'=\left(\bigcup_{i=1}^n L_ib_i \right) \cap \left(\bigcup_{i=1}^n L_i'b_i\right) =\bigcup_{i=1}^n (L_i\cap L_i')b_i.$$
Since $\CF(H)$ is closed under intersection, then $L_i\cap L_i'\in \CF(H)$, for all $i\in [n]$ and, by Proposition \ref{cf finite extension}, $K\cap K'\in \CF(G)$.

Now we prove ii. Suppose that $\CF(G)$ is closed under complement. Since it is closed under union, then it must be closed under intersection. Now, let $K\in \CF(H)$. By Lemma \ref{fi sobe}, $K\in \CF(G)$ and so $G\setminus K\in \CF(G)$. Now, $H\leq_{f.i} G$ and so $H\in \Rec(G)$. By Lemma \ref{inters cf rec}, $H\setminus K=H\cap(G\setminus K)\in \CF(G)$ and Lemma \ref{cf partial} yields that $H\setminus K\in \CF(H)$.

Conversely, suppose that  $\CF(H)$  is closed under complement and let $K\in \CF(G)$. Then $K=\bigcup_{i=1}^n L_ib_i$, for some $L_i\in \CF(H)$, $i\in[n]$ and $G\setminus K=\bigcup_{i=1}^n (H\setminus L_i)b_i$. Since $\CF(H)$ is closed under complement, then $H\setminus L_i\in \CF(H)$, for all $i\in [n]$ and by Proposition \ref{cf finite extension}, $G\setminus K\in \CF(G).$
\qed\\

Since context-free languages are not closed under intersection and complement, it is not expected for the properties in the statement of Corollary \ref{c1} to hold often. Indeed, we conjecture that $\CF(G)$ is only closed under intersection (resp. complement) if $G$ is virtually cyclic. To support these conjectures, we present the following proposition.

\begin{proposition}
Let $G$ be a finitely generated group. If $G$ is virtually abelian or virtually free, then $\CF(G)$ is closed under intersection (resp. complement) if and only if $G$ is virtually cyclic. 
\end{proposition}
\noindent\textit{Proof.} If $G$ is finite, then $\CF(G)$ is obviously closed under intersection and complement. If $G$ is virtually $\Z$, then, by \cite[Theorem 3.1]{[Her91]}, $\CF(G)=\Rat(G)$ and $\Rat(G)$ is closed under intersection and complement since $\Rat(\Z)$ is closed under intersection and complement (\cite[Propositions 3.6 and 3.9]{[Sil02b]}).

Now we will prove that if $\CF(\Z^m)$ and $\CF(F_m)$ are closed under intersection (resp. complement), then $m\leq 1$ and that suffices by Corollary \ref{c1} and the fact that every finitely generated virtually abelian group  has a free-abelian subgroup of finite index. 
%and every finitely generated virtually free group has a free subgroup of finite index. 

We start with the free-abelian case. Suppose that $G=\Z^m$, for some $m>1$. Let $A=\{e_1,\cdots, e_m\}$ be the canonical set of generators for $\Z^m$ and $\pi:\tilde A^*\to \Z^m$ be the standard surjective homomorphism. Let $K_1$ be the subset of $\Z^m$ given by the elements which have $0$ in the first component and $K_2$ be the subset of elements having $0$ in the second component.
 Then $K_1\cap K_2$ is the set of elements that have $0$ in both the first and second components. Then $$K_1\pi^{-1}=\{w\in \tilde A^*\mid n_{e_1}(w)=n_{e_1^{-1}}(w)\},$$ $$K_2\pi^{-1}=\{w\in \tilde A^*\mid n_{e_2}(w)=n_{e_2^{-1}}(w)\}$$ and $$(K_1\cap K_2)\pi^{-1}=\{w\in \tilde A^*\mid n_{e_1}(w)=n_{e_1^{-1}}(w) \text{ and }n_{e_2}(w)=n_{e_2^{-1}}(w)\}.$$
It is well known that $K_1\pi^{-1}$ and $K_2\pi^{-1}$ are context-free languages of $\tilde A^*$. Suppose that   $(K_1\cap K_2)\pi^{-1}$ is context-free, let $p$ be the constant given by the pumping lemma for context-free languages and let $z=e_1^pe_2^{p}e_1^{-p}e_2^{-p}$. Any factorization of $z$ of the form $z=uvwxy$ with $|vwx|\leq p$ and $|vx|\geq 1$ is such that $vx$ doesn't have occurrences neither of both $e_1$ and $e_1^{-1}$ nor of both $e_2$ and $e_2^{-1}$, and so $uv^2wx^2y\not\in (K_1\cap K_2)\pi^{-1}$, which contradicts the assumption that $(K_1\cap K_2)\pi^{-1}$ is context-free. Hence, $\CF(\Z^m)$ is not closed under intersection if $m>1$.

Now, suppose that $G=F_m$, for $m>1$. Let $A=\{a_1,\cdots, a_m\}$ be a free basis for $F_m$ and $\pi:\tilde A^*\to F_m$ be a surjective homomorphism. Let $K_1=\{a_1^ma_2^ma_1^{-n}\in F_m \mid m,n\geq 0\}$ and $K_2=\{a_1^ma_2^na_1^{-n}:m,n\geq 0\}$. Then $Red(K_1\pi^{-1})$ and $Red(K_2\pi^{-1})$ are context-free languages of $\tilde A^*$ but $$Red((K_1\cap K_2)\pi^{-1})=\{a_1^na_2^na_1^{-n}\in \tilde A^* \mid n\geq 0 \}$$ can easily be seen  not to be context-free by the pumping lemma. Hence, by \cite[Lemma 4.6]{[Her91]}, $K_1,K_2\in \CF(F_m)$ but $K_1\cap K_2\not\in \CF(F_m)$.

Finally suppose that $\CF(\Z^m)$ and $\CF(F_m)$ are closed under complement. Then, since $\CF(\Z^m)$ and $\CF(F_m)$ are closed under union, they must also be closed under intersection and, by the above, it follows that $m\leq 1$.

Therefore, if $G$ is a virtually abelian or virtually free group such that $\CF(G)$ is closed under complement or intersection, then  $G$ must be virtually cyclic.
\qed\\

We will now use an example from \cite{[JKLP87]}, which also appears in \cite{[Her91]}, to prove that, unlike the case of context-free subsets \cite[Corollary 4.7]{[Her91]}, algebraic subsets of $F_n$ are not closed under intersection with rational subsets of $F_n$.
\begin{example}
\label{ex free alg intersec}
Let $F_n$ be the free group of rank $n>1$ with basis $A=\{a_1,a_2,\ldots, a_n\}$.   Let $\mc G=(\{S,T\},\tilde A, P, S)$ be a context-free grammar with the following set of productions $P$:
\begin{align*}
&S\rightarrow a_1Sa_1^{-1}\\
&S\rightarrow T\\
&T\rightarrow a_1^{-1} TTa_1\\
& T\rightarrow a_2.
\end{align*}
By definition, $L(\mc G)$ is a context-free language and $K_1=L(\mc G)\pi\in \Alg(F_n).$  The set $K_2=\{a_2^n\in F_n\mid n\in \N\}$ is rational since $K_2=(a_2^*)\pi$, but $K_1\cap K_2=\{a_2^{2^n}\in F_n\mid n\in \N\}$ (see \cite[Example 1]{[JKLP87]} or \cite[Proposition 4.8]{[Her91]}). We now prove that $K_1\cap K_2\not\in \Alg(F_n)$.

Let $p:\tilde A^*\to \N^{2n}$ be the function defined by $$p(w)=(n_{a_1}(w),n_{a_1^{-1}}(w),n_{a_2}(w),n_{a_2^{-1}}(w),\ldots, n_{a_n}(w),n_{a_n^{-1}}(w)).$$
Suppose that $K_1\cap K_2=L\pi$ for some context-free language $L$. By Parikh's Theorem \cite{[Par66]}, $p(L)=\{p(w)\mid w\in L\}$ is a semilinear set of $\N^{2n}$. Denoting by $u_i$ the $i$-th coordinate of a vector $u\in \N^{2^n}$, we have that,
\begin{align}
\label{u3u4p2}
\{u_3-u_4\mid u\in p(L)\}=\{2^n\mid n\in \N\},
\end{align}
because if $w\in \tilde A^*$ is such that $w\pi=a_2^{2^k}$, then $p(w)_3-p(w)_4=2^k$.

If $p(L)$ is semilinear, then, by definition, $p(L)=\bigcup_{i=1}^k L_i$, for some linear sets $L_i$. Let $i\in[k]$ and write $$L_i=u+\N v_1 +\cdots +\N v_m,$$
for some $v_i\in \N^{2n}$, $i\in [m]$. Let $j\in [m]$. Then $u, u+v_j, u+2v_j \in L_i\subseteq p(L)$.
Hence, there are $k_0,k_1\in\N$ such that 
$$u_3-u_4=2^{k_0} \quad \text{ and  }\quad (u+v_j)_3-(u+v_j)_4=u_3+(v_j)_3-u_4-(v_j)_4=2^{k_1},$$ thus 
\begin{align}
\label{v3v4}
(v_j)_3-(v_j)_4=2^{k_1}-2^{k_0}.
\end{align}
 But now, 
\begin{align*}
(u+2v_j)_3-(u+2v_j)_4&=u_3+2(v_j)_3-u_4-2(v_j)_4\\
&= u_3-u_4 +2((v_j)_3-(v_j)_4)\\
&=2^{k_0}+2^{k_1+1}-2^{k_0+1}\\
&=2^{k_0}(2^{k_1+1-k_0}-1).
\end{align*}
Since $u+2v_j\in p(L)$, then $2^{k_0}(2^{k_1+1-k_0}-1)$ is a power of $2$ with nonnegative exponent, which implies that $k_1=k_0$. Indeed, if $k_1+1-k_0>1$, then $2^{k_1+1-k_0}-1$ is odd and if $k_1+1-k_0<1$, then $2^{k_1+1-k_0}-1$ is nonpositive.  This, together with (\ref{v3v4}), implies that $(v_j)_3-(v_j)_4=0$.
Since $j$ is arbitrary, we have that for all $j\in [m]$, $(v_j)_3-(v_j)_4=0$. Hence, for all $x\in L_i$, $x_3-x_4=u_3-u_4=2^{k_0}$.

Therefore,
$\{u_3-u_4\mid u\in p(L)\}$ is finite, which contradicts (\ref{u3u4p2}).
\end{example}

In the algebraic case, if $G$ is a finitely generated virtually abelian group, then $\Alg(G)=\Rat(G)$ and, according to \cite{[Her91]}, it is not known if the converse holds. In this case, $\Rat(G)$ is closed under intersection and complement (\cite[Propositions 3.6 and 3.9]{[Sil02b]}), and so is $\Alg(G)$. We also do not know if the converse holds, i.e., if $\Alg(G)$ being closed under intersection (resp. complement) implies that $G$ is virtually abelian. However, in the virtually free case, we can obtain the same result as in the context-free case.

\begin{corollary}
Let $G$ be a finitely generated virtually free group. Then $\Alg(G)$ is closed under intersection (resp. complement) if and only if $G$ is virtually cyclic. 
\end{corollary}
\noindent\textit{Proof.} If $G$ is virtually cyclic, then $\Alg(G)=\Rat(G)$, which is closed under intersection and complement.

If $G$ is virtually free but neither finite nor virtually $\Z$, then $G$ has a finite index free subgroup $F$ of rank greater than $1$. By Example \ref{ex free alg intersec}, $\Alg(F)$ is not closed under intersection and so by Corollary \ref{c1}, $\Alg(G)$ is not closed under intersection. 

Since $\Alg(G)$ is closed under union, then closure under complement implies closure under intersection.
\qed\\

Combining the results above on the structure of $\CF(G)$ and $\Alg(G)$ we can also prove the equivalence of a natural decidability question on $G$ and on $H$, similar to the one proved in the rational and recognizable case in \cite[Theorem 4.8]{[Sil02b]}.
\begin{proposition}\label{alg cf decidable}
Let $G=\langle A \rangle$ be a group, $H=\langle B \rangle \leq_{f.i.}G$ and $\pi_1:\tilde A^*\to G$, $\pi_2:\tilde B^*\to H$ be the canonical surjective homomorphisms. Then the following  decision problems are recursively equivalent:
\begin{enumerate}[i.]
\item taking as input a context free  grammar which generates a language $L\subseteq \tilde A^*$ such that $\pi_1(L)=K\subseteq G$, decide whether or not $K \in \CF(G)$, i.e., if $K\pi_1^{-1}\subseteq \tilde A^*$ is a context-free language.
\item taking as input a context free  grammar which generates a language $L\subseteq \tilde B^*$ such that $\pi_2(L)=K\subseteq H$, decide whether or not $K \in \CF(H)$, i.e., if $K\pi_2^{-1}\subseteq \tilde B^*$ is a context-free language.
\end{enumerate}
\end{proposition}

\noindent\textit{Proof.} Write $G$ as a disjoint union of the form $\cup_{i=1}^n Hb_i$.
Suppose that i. holds and let $K\in \Alg(H)$. Then $K\in \Alg(G)$. Since $H\leq_{f.i.}G$ and $K\subseteq H$, then $K\in \CF(G)\Leftrightarrow K\in \CF(H)$, by Lemmas \ref{cf partial} and \ref{fi sobe}.

Conversely, suppose that ii. holds and let $K\in \Alg(G)$. Then $K=\cup_{i=1}^n L_ib_i$, for some $L_i\in \Alg(H)$. By Proposition \ref{cf finite extension}, $K\in \CF(G)$ if and only if each $L_i\in \CF(H)$, for $i\in[n]$ and that we can decide from ii.
\qed\\

The analogous decision problem  for rational languages (deciding if a rational subset is recognizable) was solved in \cite[Corollary 4.9]{[Sil02b]} for finitely generated virtually abelian groups. As far as the author knows, there are no results concerning this problem in its context-free version.
Proposition \ref{alg cf decidable} shows that proving decidability for finitely generated free groups would yield decidability for finitely generated virtually free groups and, in view of Herbst's context-free analogue of Benois' Theorem, \cite[Lemma 4.6]{[Her91]}, this is equivalent to the following problem. 

\begin{problem}
Let $n\geq 2$, $X=\{x_1,\ldots, x_n\}$ and $\pi:\tilde A^*\to F_n$ be the canonical surjective homomorphism. Given a context-free grammar generating a language $L\subseteq \tilde X^*$ as input, can we decide whether $Red(L)$ is context-free or not?
\end{problem}

%%%%%%%%%%%%%%%

 % \subsection{Algebraic subsets}
%  
%\noindent\textit{Proof.} Let $K\subseteq G$ be a subset of the form $\cup_{i=1}^{n} L_ib_i$ with $L_i\in \Alg(H)$. Then, $L_i\in \Alg(G)$ for all $i\in [n]$, and so $L_ib_i\in \Alg(G)$. Since $\Alg(G)$ is closed under union, then $K\in \Alg(G)$. 
%
%Conversely, let $K\in \Alg(G)$. Then $K$ can be written as a disjoint union $$K=K\cap G=\bigcup_{i=1}^n (Hb_i\cap K).$$
%Put $K_i=Hb_i\cap K$ and let $i\in [n]$.
%Since $H\leq_{f.i.} G$, then by Proposition \ref{rec fi}, we have that $H\in \Rec(G)$, and so, $Hb_i\in \Rec(G)$, by Lemma \ref{cfletra}. It follows from  Lemma \ref{inters cf rec} that $K_i\in \Alg(G)$ and so $L_i=K_ib_i^{-1}\in \Alg(G)$ and $L_i\subseteq H$. By Lemma \ref{fi desce}, we deduce that $L_i\in \Alg(H)$ and thus
%$$K=K\cap G=\bigcup_{i=1}^n K_i=\bigcup_{i=1}^n L_ib_i.$$
%\qed\\

  %%%%%%%%%%%%%%%%%%%%%%%%%%%%%%%%%%%%%%%%%%%%%%%%%%%%%%%%%%%%%%%%%%%%%%%%%%%%%%%%%%%%%%%%%%%%%%%%%%%
  \section{Subsets of subgroups}
  \label{sec: fatou}
  The purpose of this section is to study the \emph{kind of Fatou property} in (\ref{fatou rat}) for recognizable, context-free and algebraic sets, completing the picture on this question for these four classes of subsets of a finitely generated group.
  The context-free case will lead to a new language-theoretic characterization of virtually free groups and the algebraic case will answer a question posed in \cite{[Her91]} which was further developed in \cite{[Her92]}.
   \subsection{Recognizable subsets}
 In this section we deal with the easier case of recognizable subsets and are able to describe exactly when such a property holds.
 \begin{proposition}
 Let $G$ be a group and $H\leq G$ be a subgroup. 
 Then $$\{K\subseteq H\mid K\in \Rec(G)\}\subseteq \Rec(H),$$
but the reverse inclusion holds if and only if $H\leq_{f.i.} G$.
 \end{proposition}
\noindent\textit{Proof.} Let $K\in \Rec(G)$ be such that $K\subseteq H$. Then $K$  is a (finite) union of cosets of some finite index subgroup $F\leq_{f.i.} G$:
$$K=\bigcup_{i=1}^m Fb_i,$$
for some $m\in \N$,  $b_i\in K$ such that $Fb_i\neq Fb_j$ if $i\neq j$. Since $K\subseteq H$, then 
$$K=K\cap H=\left(\bigcup_{i=1}^m Fb_i\right)\cap H=\bigcup_{i=1}^m (Fb_i\cap H).$$
Since $Fb_i\cap H$ is either empty or a coset of $F\cap H$, then 
$$K=\bigcup_{i=1}^{m'} (F\cap H)b_i',$$
for some $m'\leq m$ and $b_i'\in K$. Since $F$ has finite index in $G$, then $F\cap H$ has finite index in $H$ and $K$ can be written as a union of cosets of a finite index subgroup of $H$. Thus, $K\in \Rec(H)$.

If $[G:H]=\infty$, then $H\in \Rec(H)$, but $H\not\in \Rec(G)$, and so $\Rec(H)\not\subseteq \{K\subseteq H\mid K\in \Rec(G)\}.$ On the other hand, if  $H\leq_{f.i.} G$, then  $\{K\subseteq H\mid K\in \Rec(G)\}= \Rec(H),$ since every $K\in \Rec(H)$ is 
the union of cosets of a finite index subgroup of $H$, and so it is the union of cosets of a finite index subgroup of $G$.
\qed\\

 \subsection{Context-free subsets}
Now we deal with context-free subsets. One of the inclusions is given by Lemma \ref{cf partial}, so we only have to worry with the reverse inclusion. We will prove that it holds if and only if the group $G$ is virtually free, making use of the Muller-Schupp theorem (see \cite{[MS83]}) and of the results of the previous section concerning finite index subgroups.

We will solve the free group case first.
\begin{proposition}
\label{free case}
Let $F$ be a finitely generated free group. Then $$\CF(H)=\{K\subseteq H\mid K\in \CF(F)\},$$ for all $H\leq_{f.g.} F$.
\end{proposition}
\noindent\textit{Proof.}
Let $H\leq_{f.g} F$. By Marshall Hall's theorem (see \cite{[LS77]}) , there is a subgroup $N\leq_{f.i.} F$ such that $N=H*H'$ for some $H'\leq F$. Let $K\in \CF(H)$. We will prove that $K\in \CF(N)$, which, by Lemma \ref{fi sobe} is enough to show that $K\in \CF(F)$. So let $A$ be a free basis for $H$ and $A'$ be a free basis for $H'$ disjoint from $A$ (so $B=A\cup A'$ is a free basis for $N$) and let $\pi:\tilde A^*\to H$ and $\rho: \tilde B^*\to N$ be the standard surjective homomorphisms. So, $\rho|_{\tilde A^*}=\pi$. Since $K\in \CF(H)$, then, by \cite[Lemma 4.6]{[Her91]}, $Red(K\pi^{-1})$ is a context-free language of $\tilde A^*$. But then $Red(K\rho^{-1})=Red(K\pi^{-1})$ is context free in $\tilde B^*$, which, again by \cite[Lemma 4.6]{[Her91]}, shows that $K\in \CF(N)$.
\qed\\

Now, we can prove the main result of this subsection which provides
another language theoretical characterization of virtually free groups (see \cite{[MS83],[GHHR07],[Ant11], [DW17]} for others) among finitely generated groups.

A group is said to be \emph{Howson} if the intersection of two finitely generated subgroups is again finitely generated. We remark that free groups are Howson \cite{[How54]} and so are finite extensions of Howson groups. Hence, finitely generated virtually free groups are Howson.

\begin{theorem}\label{vfree charact}
Let $G$ be a finitely generated group. Then $G$ is virtually free if and only if $$\CF(H)=\{K\subseteq H\mid K\in \CF(G)\},$$ for all $H\leq_{f.g.} G$.
\end{theorem}
\noindent\textit{Proof.}
Let $G=\langle A\rangle$ be a finitely generated group and $\pi:\tilde A^*\to G$ be the standard surjective homomorphism. If $G$ is not virtually free, then by the Muller-Schupp theorem, we have that $\{1\}\pi^{-1}$ is not context-free, and so $\{1\}\not\in \CF(G)$. But taking $H=\{1\}$, we obviously have that $\{1\}\in \CF(H)$, and so $\CF(H)\neq\{K\subseteq H\mid K\in \CF(G)\}$

Now, suppose that $G$ is virtually free.  Then $G$ admits a decomposition as a disjoint union of the form $$G=Fb_1\cup \cdots \cup Fb_m,$$ for some free group $F\leq_{f.i.}G$ and $b_i\in G$, for $i\in[m]$.
By Lemma \ref{cf partial}, we only have to prove that $\CF(H)\subseteq \CF(G)$. Let $H\leq_{f.g} G$ and $K\in \CF(H)$. We will show that $K\in \CF(G)$. 

We have that $$H=H\cap G= H\cap \left(\bigcup_{i=1}^m Fb_i\right)= \bigcup_{i=1}^m(H\cap Fb_i) =\bigcup_{i=1}^{m'}(H\cap F)b_i' ,$$
for some $m'\leq m$ and $b_i'\in H$, since $H\cap Fb_i$ is either empty or a coset of $H\cap F$. Moreover, we can assume the cosets $(H\cap F)b_i'$ to be disjoint. Since $K\in \CF(H)$, then 
$$K=\bigcup_{i=1}^{m'} L_i b_i',$$
for some $L_i\in \CF(H\cap F)$, $i\in[m']$. Since virtually free groups are Howson, then $H\cap F \leq_{f.g} F$ and, by Proposition \ref{free case}, it follows that $L_i\in \CF(F)$, for all $i\in [m']$. By Lemma \ref{fi sobe}, this implies that $L_i\in \CF(G)$. Using Lemma \ref{cfletra}, we get that $L_ib_i'\in \CF(G)$, for all $i\in[m']$ and since $\CF(G)$ is closed under union, then $K\in \CF(G)$. 
\qed\\

  \subsection{Algebraic subsets}
  Finally, we tackle the same problem for algebraic subsets.
   \begin{problem}
 \label{question}
 Let $G$ be a group and $H\leq G$ be a subgroup. Is it true that $\Alg(H)=\{K\subseteq H\mid K\in \Alg(G)\}$?
 \end{problem}
  Notice that, similarly to what happens in the rational case, in the algebraic case the inclusion $\Alg(H)\subseteq \{K\subseteq H\mid K\in \Alg(G)\}$ is obvious, so we will only deal with the reverse inclusion.

  This particular question was raised by Herbst in \cite{[Her91]} and answered affirmatively in the particular case where $G$ is a virtually free group in \cite{[Her92]}. We also know that the answer is positive in the case where $H$ is a finite index subgroup of $G$ by Lemma \ref{fi desce}. However,  we will prove that this is not always the case by constructing a specific counterexample.

 Let $A=\{a_1,\ldots, a_n\}$, $G=\langle A\mid R\rangle$ be a group, $\varphi\in \Aut(G)$ and $g\in G$. The orbit of $g$ through $\varphi$ is the set $\{g\varphi^k\mid k\in \N\}$ and we denote it by $\Orb_\varphi(g)$.
 
 Consider the semidirect product 
\begin{align}
\label{presentation semidirect}
G\rtimes_\varphi \Z=\langle A,t \mid R, t^{-1}a_it= a_i\varphi\rangle.
\end{align}

Using the relations, every element of $G\rtimes_\varphi \Z$ can be rewritten as an element of the form $t^ag$, where $a\in \Z$ and $g\in G$, in a unique way. 

\begin{remark}
In a group of the form $G\rtimes_\varphi \Z$, $G\pi^{-1}=\{w\in \widetilde{A\cup\{t\}}^*\mid n_t(w)=n_{t^{-1}}(w)\}$ is context-free, and so $G\in \CF(G\rtimes_\varphi \Z)$. In particular, since $G$ has infinite index in $G\rtimes_\varphi \Z$, this is an example of a group where recognizable and context-free subsets do not coincide.
\end{remark}
 
 \begin{theorem}
 \label{thm}
Let $G$ be a finitely generated group  and $\varphi\in \Aut(G)$. If $\Alg(G)=\{K\subseteq G\mid K\in \Alg(G\rtimes_\varphi \Z)\}$, then the orbit through $\varphi$ of every element is an algebraic subset of $G$.
 \end{theorem}
 \noindent\textit{Proof.}  Let $G$ be a group generated by a finite set $A=\{a_1,\ldots a_n\}$ and $\pi:\tilde A^*\to G$ be the standard surjective homomorphism. Then $G\rtimes_\varphi \Z$ admits a presentation of the form (\ref{presentation semidirect}) and there is a natural surjective homomorphism $\rho:\widetilde{A\cup \{t\}}^*\to G\rtimes_\varphi\Z$ such that $\rho|_A=\pi$ (identifying $G$ with the subset $\{t^0g\mid g\in G\}$).  For every $w\in \tilde A^*$, the language $L=\{t^{-n}wt^n\mid n\in \N\}\subseteq \widetilde{A\cup \{t\}}^*$ is context-free and so $L\rho$ is an algebraic subset of $G\rtimes_\varphi \Z$. But $(t^{-n}wt^{n})\rho=t^0(w\pi)\varphi^n$, for all $n\in \N$, thus $L\rho=\{t^0(w\pi)\varphi^n\mid n\in \N\}\subseteq G$.  Since $\Alg(G)=\{K\subseteq G\mid K\in \Alg(G\rtimes_\varphi \Z)\}$, then $\Orb_\varphi(w\pi)\in \Alg(G)$, for every $w\in \tilde A^*$. 
 \qed\\
 
 This allows us to construct a counterexample to Problem \ref{question}.

 \begin{example}
 Consider the group $\Z^2$ and $ Q=
\begin{bmatrix}
    2 &1 \\
    1& 1\end{bmatrix}
\in \text{GL}_m(\Z)$ 

Since $\Z^2$ is abelian, we have that $\Alg(\Z^2)=\Rat(\Z^2)$. We will see that the orbit of $(1,0)\in \Z^2$ is not rational, which by Theorem \ref{thm}, implies that $\Alg(\Z^2)\neq \{K\subseteq G\mid K\in \Alg(\Z^2\rtimes_Q \Z)\}$.
Let $(f_n)_n\in \N$ be the Fibonacci sequence.

We have that, for $k\in \N$
\begin{align*}
Q^k=\begin{bmatrix}
    f_{2k+1} &f_{2k} \\
    f_{2k}& f_{2k-1}\end{bmatrix}
 \end{align*}
 This can be seen by induction on $k$. It holds for $k=1$. Suppose that it holds for all integers up to some $r\in \N$. Then 
 $$Q^{r+1}=Q^rQ=\begin{bmatrix}
    f_{2r+1} &f_{2r} \\
    f_{2r}& f_{2r-1}\end{bmatrix}\begin{bmatrix}
    2 &1 \\
    1& 1\end{bmatrix}=\begin{bmatrix}
    2f_{2r+1}+f_{2r} &2f_{2r}+f_{2r-1} \\
    f_{2r}+f_{2r+1}& f_{2r}+f_{2r-1}\end{bmatrix}=\begin{bmatrix}
    f_{2r+3} &f_{2r+2} \\
    f_{2r+2}& f_{2r+1}\end{bmatrix}.$$
The orbit of $(1,0)$ is the set of the first rows of powers of $Q$:
$$\Orb_Q((1,0))=\{(1,0)\}\cup\{(f_{2k+1},f_{2k})\mid k\in \N\}
.$$

Suppose that $\Orb_Q((1,0))$ is rational. Then so is its projection to the first component $\{f_{2k+1}\mid k\in \N\}$. Let $A=\{a,a^{-1}\}$. By Benois' Theorem, the language $L=\{a^{f_{2k+1}}\mid k\in \N\}\subseteq  A^*$ is rational, which is absurd by the pumping lemma for rational languages.
 \end{example}

In the same spirit, we can see that the Baumslag-Solitar group $BS(2,1)$ is also a counterexample to Problem \ref{question}.

 \begin{example}\label{ex baumslag}
Consider the Baumslag-Solitar group $G=BS(2,1)=\langle a,t \mid ta^2t^{-1}=a\rangle.$ Let $H=\langle a\rangle$. Then $S=\{t^{-n}at^n\}=\{a^{2^n}\}\subseteq H$ and $S\in \Alg(G)$, but $S\not\in \Alg(H)$, since $H\simeq \Z$ and $S\not\in \Rat(H)=\Alg(H)$.
 \end{example}
 
\begin{remark}
We know from \cite{[Her92]} that if $G$ is virtually free, then $\Alg(H)= \{K\subseteq H\mid K\in \Alg(G)\}$, for every $H\leq_{f.g.} G$. However, unlike the context-free case, we cannot expect this property to characterize virtually free groups, since, for example in abelian groups, the set of algebraic subsets and rational subsets coincide and property (\ref{fatou rat}) always holds.
\end{remark}

\section*{Acknowledgements}
The author is grateful to Pedro Silva for fruitful discussions of these topics, which  improved the paper, and to the anonymous referee for several corrections and for Example \ref{ex baumslag}.
The author was  supported by the grant SFRH/BD/145313/2019 funded by Funda\c c\~ao para a Ci\^encia e a Tecnologia (FCT).%
\bibliographystyle{plain}
\bibliography{Bibliografia}
\end{document}